\documentclass[onefignum]{siamart251216}

\pdfoutput=1
\usepackage{nmath}
\usepackage{xparse}
\usepackage{algorithm}
\usepackage{algpseudocodex}
\usepackage[shortlabels]{enumitem}

\usepackage{graphicx}
\graphicspath{{fig/}}

\def\Figure#1#2{%
\centering%
\includegraphics{#1}%
\caption{#2}%
\label{fig:#1}%
}

\usepackage{float} \floatstyle{plain} \newfloat{multalg}{ht}{multalg}

\def\algsize{\footnotesize}
\def\algspace{1.2}
\usepackage{xspace}
\def\explicit{\texttt{explicit}\xspace}
\def\hybridCS{\texttt{hybridCS}\xspace}
\def\hybridRCM{\texttt{hybridRCM}\xspace}
\def\explicitRCM{\texttt{explicitRCM}\xspace}

\def\deq{\doteq}

\def\h{e} 

\def\P{{\varphi}} 

\NewDocumentCommand{\seq}{m O{#3} m m m}{%
  \left\{#1_{#2}\right\}_{#3=#4}^{#5}
}

\newcommand{\ze}[2]{{\zeta_{#1}^{#2}}}

\usepackage{hyperref}

\usepackage[nameinlink]{cleveref}
\crefformat{equation}{(#2#1#3)}

\crefname{appendix}{appendix}{appendices}
\Crefname{appendix}{Appendix}{Appendices} 

\crefname{algorithm}{algorithm}{algorithms}
\Crefname{algorithm}{Algorithm}{Algorithms}

\makeatletter
\AddToHook{cmd/appendix/before}{\crefalias{section}{appendix}}\makeatother

\algrenewcommand\algorithmicrequire{\textbf{Input:}}
\def\Input{\Require}

\def\fftwpp{{\tt FFTW++}~}

\def\fft#1{{\tt fft}(#1)}
\def\rcfft#1{{\tt rcfft}(#1)} 

\def\ifft#1{{\tt ifft}(#1)}

\def\ircfft#1{{\tt crfft}(#1)}

\def\convolveRCMone{\texttt{convolveRCM1}}
\def\convolveRCMtwo{\texttt{convolveRCM2}}

\def\zeropad{\texttt{zeroPad}}
\def\complexpack{\texttt{complexPack}}
\def\complexunpack{\texttt{complexUnPack}}
\def\fftx{\texttt{fftx}}
\def\ifftx{\texttt{ifftx}}
\def\ffty{\texttt{ffty}}
\def\iffty{\texttt{iffty}}

\def\True{\textbf{True}}
\def\False{\textbf{False}}

\def\Or{\textbf{or}}

\def\Not{\textbf{not}}

\def\Forward{{\tt Forward}}
\def\Backward{{\tt Backward}}

\def\forwardReal{{\tt forward1R}}

\def\backwardReal{{\tt backward1R}}

\def\param{L,m,q}

\def\ForwardIs#1{$\Forward\gets~$\hyperlink{alg:#1}{\csname #1\endcsname}}
\def\BackwardIs#1{$\Backward\gets~$\hyperlink{alg:#1}{\csname #1\endcsname}}

\def\benchmarkCaption#1#2#3{In-place #1 real convolutions of two arrays of #2 on #3.}
\def\benchCaptionOne#1{\benchmarkCaption{1D}{length~$L$}{#1}}
\def\benchCaptionTwo#1{\benchmarkCaption{2D}{size~$L\times L$}{#1}}
\def\benchCaptionThree#1{\benchmarkCaption{3D}{size~$L\times L\times L$}{#1}}

\def\first{\texttt{first}}
\def\multBinaryRCM{\texttt{multRCM}}

\def\FBdesc#1#2#3#4{is the #4 #1 transform for residue #2 when #3.}

\def\Fdesc#1#2#3{\FBdesc{forward}{#1}{\hbox{#2}}{#3}}

\def\Bdesc#1#2#3{\FBdesc{backward}{#1}{\hbox{#2}}{#3}}

\def\AA{W}

\usepackage{xcolor}

\let\Re\relax
\DeclareMathOperator{\Re}{Re}
\let\Im\relax
\DeclareMathOperator{\Im}{Im}

\usepackage[
  style=ext-numeric,
  doi=false,
  url=false,
  isbn=false,
  giveninits=true,
  articlein=false,
  sorting=nyt,
  sortcites
]{biblatex}
\DeclareFieldFormat{pages}{#1}

\headers{Hybrid Dealiasing for Real Convolutions}{Noel Murasko and John C.\ Bowman}

\title{Hybrid Dealiasing and Implicit Packing for Real Convolutions}
\author{Noel Murasko\thanks{University of Alberta
(\email{murasko@ualberta.ca}, \email{bowman@ualberta.ca})} \and John C.\ Bowman\footnotemark[1]}

\date{}

\begin{document}
\maketitle
\noindent
\begin{abstract}
Hybrid dealiasing is an FFT-based method for computing linear convolutions of
complex-valued data that reduces the cost of dealiasing by performing
zero padding implicitly. We develop two new algorithms that extend hybrid
dealiasing to real-valued convolutions.

The first algorithm exploits conjugate symmetries in the transformed data and
computes each residue contribution directly. The second algorithm employs
complex-valued hybrid dealiasing via a new \emph{implicit packing} technique,
which packs real-valued data into complex-valued arrays and performs all
unpacking and packing operations implicitly in the transformed domain.
Multidimensional versions of both algorithms are obtained by recursive
decomposition into lower-dimensional convolutions.

Based on these algorithms, efficient routines to compute convolutions of
equal-length input arrays are implemented in the open-source {\tt FFTW++}
library. These routines outperform the standard method of explicit
zero padding in one, two, and three dimensions. In many cases,
implicit packing without hybrid dealiasing is
nearly as fast as the general hybrid dealiasing algorithm,
particularly when multithreading.
\end{abstract}

\begin{keywords}
dealiasing, hybrid dealiasing, implicit packing,
real-valued convolution, discrete Fourier transform, implicit dealiasing,
zero padding, fast Fourier transform
\end{keywords}

\begin{AMS}
65T50, 65Y05
\end{AMS}

\section{Introduction}\label{sec:introduction}
The fast Fourier transform (FFT) can be used to compute discrete circular
convolutions efficiently; however, many applications require linear
convolutions. A naive application of FFTs to linear convolutions produces data
polluted by aliasing errors arising from the assumed periodicity. The process
of removing these aliases is called \emph{dealiasing}.

The standard dealiasing method explicitly zero pads the input arrays; we call
this \emph{explicit dealiasing}. While effective, explicit dealiasing
requires reading and multiplying values known \textit{a priori} to be zero.
This inefficiency worsens with dimension: explicitly dealiasing a~$d$-dimensional convolution of two arrays of equal size requires approximately~${(2^{d}-1)}/2^d$ of the inputs to be zero.

An alternative is \emph{implicit dealiasing} \cite{BowmanRoberts2011,RobertsBowman2018}, which constructs padded and unpadded FFTs to account for known zeros implicitly, eliminating the need for explicit padding. Implicit dealiasing often outperforms explicit dealiasing, with the greatest gains in higher dimensions. However, the original formulation was limited, as it did not handle arbitrary padding ratios (the ratio of output FFT size to input size) or convolutions of real-valued data.

\emph{Hybrid dealiasing} improves upon implicit dealiasing by combining explicit and
implicit zero padding. In \cite{MuraskoBowman2024}, we developed hybrid
dealiasing algorithms for complex and Hermitian-symmetric arrays, generalizing
implicit dealiasing to arbitrary padding ratios. These algorithms were
implemented in version 3 of \fftwpp \cite{fftwpp} and were shown to match or
outperform implicit dealiasing in one, two, and three dimensions. A summary is
given in \cref{sec:recap}.

Here we extend hybrid dealiasing to efficiently compute convolutions of real-valued data via two
new algorithms. The first (\cref{sec:symmetries}) is analogous to the
Hermitian-symmetric algorithm of \cite{MuraskoBowman2024}: the convolution is
decomposed into \emph{residues} (see \cref{sub:oneRes}), and each residue
contribution is computed by exploiting conjugate symmetries in the Fourier
transformed data.

The second algorithm (\cref{sec:rcm}) packs real-valued data into complex arrays
of half the size, enabling use of the efficient complex hybrid dealiasing
algorithms of \cite{MuraskoBowman2024}. Standard packing methods require an
expensive unpacking and packing step in Fourier space. We avoid this via
\emph{implicit packing}, which performs these operations implicitly through the
multiplication operator.

In \cref{sec:multidim}, we extend both algorithms to multidimensional
convolutions using the decomposition of \cite{MuraskoBowman2024}. This
decomposition computes a multidimensional convolution recursively, yielding
substantial memory savings for large arrays. Extending the implicit packing
algorithm is particularly involved, as the modified multiplication operator
is not an element-wise operation.

Both algorithms are implemented in version 3 of \fftwpp in one, two, and three
dimensions. We benchmark them in \cref{sec:results} against explicitly
dealiased convolutions using real-to-complex/complex-to-real FFTs. Both
hybrid-dealiased algorithms outperform explicit dealiasing, with implicit
packing typically being fastest. We also compare our algorithms to a routine using
implicit packing with explicit dealiasing, which is competitive in some cases.

In \cref{sec:conclusion} we discuss potential future work, including application
to convolutional neural networks.

\section{A recap of complex hybrid dealiasing}\label{sec:recap}

In this section, we summarize the basic ideas of hybrid dealiasing in one dimension (see \cite[\S
2]{MuraskoBowman2024} for details). Given a one-dimensional complex array
$\vf=\seq{f}{j}{0}{L-1}$ that must be zero padded to at least length~$M\geq L$ before
computing its DFT $\vF\deq \seq{F}{k}{0}{M}$ (we denote DFTs with
upper-case letters and definitions by $\deq$), we choose a parameter $m\in \N$ and set $p\deq \ceil{L/m}$ and~$q\deq \ceil{M/m}$. We explicitly zero pad~$\vf$ from length~$L$ to~$pm$, then compute the implicitly zero-padded DFT of length~$qm$.

More specifically, we reindex $\vf$ and~$\vF$ using
\begin{equation*}
j=tm+s,\quad t=0,\ldots, p-1, \quad s=0,\ldots, m-1,
\end{equation*}
and
\begin{equation*}
k=q\ell+r,\quad \ell=0,\ldots, m-1, \quad r=0,\ldots, q-1.
\end{equation*}
For~$N\in \N$, let~$\zeta_N\deq e^{2\pi i/N}$. Then, following the Cooley--Tukey
decomposition \cite{CooleyTukey1965}, the implicitly zero-padded DFT of~$\vf$ is
given by
\begin{equation}\label{eqn:forwardr}
F_{q\ell+r} =\sum_{s=0}^{m-1}\zeta_{m}^{\ell s}\zeta_{qm}^{rs}\sum_{t=0}^{p-1}\zeta_{q}^{rt}f_{tm+s},
\end{equation}
and the inverse is given by
\begin{equation}\label{eqn:backwardr}
f_{tm+s}=\frac{1}{qm}\sum_{r=0}^{q-1}\zeta_{q}^{-tr}\zeta_{qm}^{-sr}\sum_{\ell=0}^{m-1}\zeta_{m}^{-s\ell}F_{q\ell+r}.
\end{equation}
These transforms require $q$ FFTs of length~$m$. We note that the parameter $m$ can be optimized to minimize the total convolution cost.

\subsection{One residue at a time}\label{sub:oneRes}
For each \emph{residue}~$r\in\{0,\ldots,q-1\}$, define the \emph{residue contribution}
$\vF_r\doteq\seq{F}[q\ell + r]{\ell}{0}{m-1}$. When unambiguous, we use \emph{residue} as shorthand for \emph{residue contribution}; for example, computing residue~$r$ means computing~$\vF_r$.

Hybrid-dealiased convolutions can be computed one residue at a time
\cite[\S 2.2]{MuraskoBowman2024}. For complex-valued data this is primarily a
memory-saving optimization; for the real algorithm of \cref{sec:symmetries} it
is essential, since not all residue contributions are computed identically.

From \cref{eqn:forwardr}, $\vF_r$ is computed by applying an FFT of length~$m$ to
$$
\zeta_{qm}^{rs}\sum_{t=0}^{p-1}\zeta_{q}^{rt}f_{tm+s}.
$$
For the inverse, we apply another FFT of length~$m$ to obtain
$\vh_{r}=\{h_{r,j}\}_{j=0}^{pm-1}$ where
\begin{equation}\label{eqn:accumulator}
h_{r,tm+s}\doteq
\ze{q}{-tr}\ze{qm}{-sr}\sum_{\ell=0}^{m-1}\ze{m}{-s\ell}F_{q\ell+r}.
\end{equation}
Following \cref{eqn:backwardr}, we accumulate over residues to obtain the inverse
\begin{equation}\label{eqn:accumulation}
\vf = \frac{1}{qm}\sum_{r=0}^{q-1}\vh_r.
\end{equation}

Each residue contribution to the implicitly-padded complex DFTs in
\cref{eqn:forwardr,eqn:backwardr} requires an FFT of length~$m$, which can be computed by an existing library such as \texttt{FFTW} \cite{fftw,FrigoJohnson2005}.

\subsection{Hermitian symmetric and real arrays}\label{sub:Hermitian}
A periodic array~$\vf=\seq{f}{j}{0}{L-1}$ with period~$L$ is Hermitian symmetric if $f_j = \conj{f_{L-j}}$ (where the bar denotes complex conjugation) for all~$j\in \{0,\ldots, L-1\}$. This symmetry implies that the \textit{DC mode}~$f_0$ is real; if $L$ is even, the \textit{Nyquist mode}~$f_{L/2}$ is also real.

It is easy to show that an array is Hermitian symmetric if and only if its DFT is real. Since real-valued data appears in many applications, considerable effort has been devoted to efficient real-to-complex and complex-to-real FFTs.

As in the complex case, hybrid-dealiased Hermitian-symmetric convolutions can be computed one residue at a time using complex-to-real/real-to-complex FFTs of length~$m$ \cite[\S 4]{MuraskoBowman2024}. The argument is straightforward: if $\vf$ is Hermitian symmetric, its DFT~$\vF=\seq{F}{k}{0}{qm-1}$ is real. For any~$r\in\{0,\ldots, q-1\}$, the residue contribution~$\vF_r=\seq{F}[q\ell+r]{\ell}{0}{m-1}$ is also real (being a sub-array of~$\vF$), so $\vF_r$ must be the DFT of a Hermitian symmetric array.

Unfortunately, this approach does not extend to implicitly-padded real transforms. If~$\vg=\seq{g}{j}{0}{qm-1}$ is real, then~$\vG=\seq{G}{k}{0}{qm-1}$ is Hermitian symmetric. For any~$r\in\{1,\ldots, q-1\}$, the residue contribution~$\seq{G}[q\ell+r]{\ell}{0}{m-1}$ is not generally Hermitian symmetric (at $\ell=0$, the DC mode~$\vG_{r}$ need not be real), so it cannot be the DFT of a real array. Consequently, real hybrid dealiasing cannot rely on real-to-complex/complex-to-real FFTs.

\section{Real hybrid dealiasing via conjugate symmetries}\label{sec:symmetries}
In this section, we assume that $\vf=\seq{f}{j}{0}{pm-1}$ is real-valued. Unless stated otherwise, the notation and indices of \cref{sec:recap} apply throughout this section.

\subsection{Exploiting conjugate symmetries}\label{sub:conjugate}
Although real-to-complex/complex-to-real FFTs cannot be used directly, we can exploit conjugate symmetries in the data. By Hermitian symmetry, the complex conjugate of~$F_{q\ell+r}$ satisfies
$$\conj{F_{q\ell+r}}
=F_{qm-(q\ell+r)}=F_{q(m-\ell)-r},$$
so \cref{eqn:accumulator} satisfies
\begin{align}\label{eqn:hrHerm}
\begin{split}
\conj{h_{r,tm+s}}&=\zeta_{q}^{tr}\zeta_{qm}^{sr}\sum_{\ell=0}^{m-1}\zeta_{m}^{s\ell} \conj{F_{q\ell+r}}=\zeta_{q}^{tr} \zeta_{qm}^{sr}\sum_{\ell=0}^{m-1}\zeta_{m}^{s\ell} F_{q(m-\ell)-r}\\
&= \zeta_{q}^{-t(q-r)}\zeta_{qm}^{-s(q-r)}\sum_{\ell=0}^{m-1}\zeta_{m}^{-s\ell} F_{q\ell+q-r}=h_{q-r,tm+s},
\end{split}
\end{align}
where the penultimate equality follows from $\ell\mapsto m-1-\ell$. This yields
\begin{equation}\label{eqn:real-accumulation}
\vf=\frac{1}{qm}\left( \vh_{0}+2\sum_{r=1}^{\ceil{q/2}-1}\Re \vh_{r}+\vh_{q/2} \right),
\end{equation}
where we take~$\vh_{q/2}\deq \vector{0}$ when~$q$ is odd.
\Cref{eqn:real-accumulation} forms the basis of the real hybrid dealiasing algorithm.
There are three cases to consider.

\subsubsection{Case \texorpdfstring{$\vh_0$}{r=0}}\label{ssub:zero}
This case is straightforward: since $\sum_{t=0}^{p-1}f_{tm+s}$ is real, the forward transform
\begin{equation*}
F_{q\ell} =\sum_{s=0}^{m-1}\zeta_{m}^{\ell s} \sum_{t=0}^{p-1}f_{tm+s}
\end{equation*}
can be computed with a real-to-complex FFT of length~$m$. Similarly, \cref{eqn:accumulator} becomes
\begin{equation*}
h_{0,tm+s}= \sum_{\ell=0}^{m-1}\zeta_{m}^{-s\ell}F_{q\ell},
\end{equation*}
which requires a complex-to-real FFT of length~$m$.

\subsubsection{Case \texorpdfstring{$\vh_r$, $r\in\{1,\ldots,\ceil{q/2}-1\}$}{r=1,...,ceil(q/2)-1}}
This case is also straightforward; since \cref{eqn:real-accumulation}
does not require $\vh_{q-r}$, this formulation uses complex
FFTs of length~$m$ without sacrificing efficiency.

\subsubsection{Case \texorpdfstring{$\vh_{q/2}$}{r=q/2}}\label{ssub:q2}
This case requires more care and arises only when $q$ is even. Simplifying \cref{eqn:forwardr}, the residue~$\vF_{q/2}$ is given by
\begin{equation}\label{eqn:resq2}
F_{q\ell+q/2} =\sum_{s=0}^{m-1}\zeta_{m}^{\ell s} \zeta_{2m}^{s}\sum_{t=0}^{p-1}(-1)^tf_{tm+s}.
\end{equation}
Since $\zeta_{2m}^{s}\sum_{t=0}^{p-1}(-1)^tf_{tm+s}$ is not real, a real-to-complex FFT (as used for $r=0$) is inapplicable. A full complex FFT of length~$m$ would also be wasteful.

The key observation is that $\vF_{q/2}$ is \emph{almost}\footnote{True Hermitian symmetry requires
$\conj{F_{q\ell+q/2}}=F_{q(m-\ell)+q/2}$.} Hermitian
symmetric: the $q/2$ residue contribution is closed under conjugation:
\begin{equation}\label{eqn:realq2symmetry}
\conj{F_{q\ell+q/2}} = F_{qm-(q\ell+q/2)}=F_{q\left(m-1-\ell\right)+q/2}.
\end{equation}
That is, the conjugates of the even indices of~$\vF_{q/2}$ coincide with its
odd indices.
We restrict $m$ to be even, let $\h\deq m/2$, and reindex $s$:
$$s=a \h +b, \quad a\in\{0,1\}, \ b\in\{0, \ldots, \h-1\}.$$
Only the even\footnote{The choice of even indices is arbitrary; the odd indices yield identical results.}
indices of~$\vF_{q/2}$ need to be computed. For $c=0,\ldots, \h-1$, the even indices of \cref{eqn:resq2} are
\begin{equation}\label{eqn:forwardq2}
F_{2cq+q/2}=\sum_{b=0}^{\h-1}\zeta_{\h}^{cb} \zeta_{2m}^{b}\sum_{t=0}^{p-1}(-1)^t\left(f_{tm+b}+if_{tm+\h+b}\right).
\end{equation}
To compute the inverse, define
\begin{equation}\label{eqn:wb}
w_{b}\deq \zeta_{2m}^{-b}\sum_{c=0}^{\h-1}\zeta_{\h}^{-bc}F_{2cq+q/2}.
\end{equation}
Ignoring the twiddle factors, $\vw$ is essentially the inverse DFT of
\cref{eqn:forwardq2}, and can be computed with an FFT of length~$\h$. To recover~$\vh_{q/2}$ from~$\vw$, apply \cref{eqn:realq2symmetry} to the conjugate of \cref{eqn:wb}:
\begin{align*}\conj{w_{b}}&=\zeta_{2m}^{b}\sum_{c=0}^{\h-1}\zeta_{\h}^{bc}\conj{F_{2cq+q/2}}=\zeta_{2m}^{b}\sum_{c=0}^{\h-1}\zeta_{\h}^{bc}F_{(2\h - 1 - 2c)q+q/2}\\
&=\zeta_{2m}^{b}\sum_{c=0}^{\h-1}\zeta_{\h}^{-b(c+1)}F_{(2c +1)q+q/2}=\zeta_{2m}^{-3b}\sum_{c=0}^{\h-1}\zeta_{\h}^{-bc}F_{(2c +1)q+q/2},
\end{align*}
where the penultimate equality follows from $c\mapsto \h-1-c$. Then \cref{eqn:accumulator} becomes
\begin{align*}
\begin{split}
h_{q/2,tm+a\h +b}&=(-1)^t\zeta_{2m}^{-(ae+b)}\sum_{\ell=0}^{m-1}\zeta_{m}^{-(ae+b)\ell}F_{q\ell+q/2}\\
&=(-1)^ti^{-a}\zeta_{2m}^{-b}\left[\sum_{c=0}^{\h-1}\zeta_{\h}^{-bc}F_{2cq+q/2}+(-1)^a\zeta_{m}^{-b}\sum_{c=0}^{\h-1}\zeta_{\h}^{-bc}F_{(2c +1)q+q/2}\right]\\
&=(-1)^ti^{-a}\left[w_{b}+(-1)^a \conj{w_{b}}\right].
\end{split}
\end{align*}
Considering~$a=0$ and~$a=1$ separately gives us
\begin{equation}\label{eqn:accumulateq2}
h_{q/2,tm+b}=(-1)^t2\Re w_{b}, \quad h_{q/2,tm+\h+b}=(-1)^t2\Im w_{b},
\end{equation}
and so we can compute~$\vh_{q/2}$ using a complex FFT of length~$\h$.

\subsection{Summary}
The algorithm described in this section accumulates residue contributions via \cref{eqn:real-accumulation} as follows:
\begin{itemize}
  \item When~$r=0$, use a real-to-complex/complex-to-real FFT of length~$m$.
  \item When~$r\in \{1,\ldots, \ceil{q/2}-1\}$, use a complex FFT of length~$m$,
  simultaneously yielding the contributions for residues~$r$ and~$q-r$.
  \item When~$r=q/2$, require that~$m$ be even and use a complex FFT of length~$m/2$.
\end{itemize}
Pseudocode for the forward and backward transforms is given in
\cref{alg:forwardReal,alg:backwardReal} for the case $p=1$. The case $p=2$ appears in \cite[Algorithms 14 and 15]{Murasko2024_thesis}. For $p>2$, we apply an optimization that computes the preprocessing and postprocessing using FFTs of length~$p$, analogous to \cite[\S6.1]{MuraskoBowman2024}; see \cref{sec:realInner} for details.

\begin{multalg}[htbp]
\begin{minipage}[t]{0.49\textwidth}
\centering
\begin{algorithm}[H]
\algsize
\caption{\forwardReal~\Fdesc{$r$}{$p=1$}{real}}
\begin{algorithmic}\label{alg:forwardReal}
\setlength{\baselineskip}{\algspace\baselineskip}
\Input{$\seq{f}{j}{0}{L-1},\param,r$}
\State{$\h\gets \floor{m/2}$}
\If{$r<q/2$}
  \For{$s=0,\ldots, L-1$}
    \State{$W_s\gets \ze{qm}{rs}f_s$}
  \EndFor
  \For{$s=L,\ldots, m-1$}
    \State{$W_s\gets 0$}
  \EndFor
  \If{$r=0$}
    \State{$\seq{V}{\ell}{0}{\h+1} \gets \rcfft{\seq{W}{s}{0}{m-1}}$}
    \Return{$\seq{V}{k}{0}{\h+1}$}
  \Else
    \State{$\seq{V}{\ell}{0}{m-1} \gets \fft{\seq{W}{s}{0}{m-1}}$}
    \Return{$\seq{V}{k}{0}{m-1}$}
  \EndIf
\ElsIf{$r=q/2$}
  \State{$B_1\gets \max(0,L-\h)$}
  \State{$B_2\gets \min(\h,L)$}
  \For{$b = 0, \ldots, B_1-1$}
    \State{$W_b\gets \ze{2m}{b}\left(f_b + if_{\h+b}\right)$}
  \EndFor
  \For{$b =B_1, \ldots, B_2-1$}
    \State{$W_b\gets \ze{2m}{b}f_b$}
  \EndFor
  \For{$b = L, \ldots, \h -1$}
    \State{$W_b\gets 0$}
  \EndFor
  \State{$\seq{V}{c}{0}{\h-1} \gets \fft{\seq{W}{b}{0}{\h-1}}$}\\
  \Return{$\seq{V}{k}{0}{\h-1}$}
\EndIf
\end{algorithmic}
\end{algorithm}
\end{minipage}
\hfill
\begin{minipage}[t]{0.49\textwidth}
\centering
\begin{algorithm}[H]
\algsize
\caption{\backwardReal~\Bdesc{$r$}{$p=1$}{real} The normalization factor $1/(qm)$ from \cref{eqn:real-accumulation} would be applied during the final accumulation step.}
\begin{algorithmic}\label{alg:backwardReal}
\setlength{\baselineskip}{\algspace\baselineskip}
\Input{$\seq{F}{k}{0}{2m-1},\param, r$}
\State{$\h\gets \floor{m/2}$}
\If{$r=0$}
  \State{$\seq{W}{s}{0}{m-1} \gets \ircfft{\seq{F}{\ell}{0}{\h+1}}$}
\ElsIf{$r<\ceil{q/2}$}
  \State{$\seq{W}{s}{0}{m-1} \gets \ifft{\seq{F}{\ell}{0}{m-1}}$}
  \For{$s=0,\ldots, L-1$}
    \State{$W_s\gets 2\Re{\ze{qm}{-rs}W_s}$}
  \EndFor
\Else
   \State{$\seq{W}{b}{0}{\h-1} \gets \ifft{\seq{F}{c}{0}{\h-1}}$}
  \State{$B_1\gets \max(0,L-\h)$}
  \State{$B_2\gets \min(\h,L)$}
  \For{$b = 0, \ldots, B_1-1$}
    \State{$W_{\h+b}\gets 2\Im{\ze{2m}{-b}W_{b}}$}
    \State{$W_{b}\gets 2\Re{\ze{2m}{-b}W_{b}}$}
  \EndFor
  \For{$b = B_1, \ldots, B_2-1$}
    \State{$W_{b}\gets 2\Re{\ze{2m}{-b}W_{b}}$}
  \EndFor
\EndIf
\Return{$\seq{W}{s}{0}{L-1}$}
\end{algorithmic}
\end{algorithm}
\end{minipage}
\end{multalg}

\section{Real hybrid dealiasing via implicit packing}\label{sec:rcm}
Our second algorithm packs real-valued data into a complex array of half the size. Let~$\Nh\in\N$. For a real array~$\vf\deq\seq{f}{j}{0}{2\Nh-1}$, define its \emph{complex packing} $\vtf=\{\tf_j\}_{j=0}^{\Nh-1}$ by
\begin{equation}\label{eqn:complexPacking}
\tf_j=f_{2j}+if_{2j+1}.
\end{equation}
That is, the even and odd indices of~$\vf$ form the real and imaginary parts of~$\vtf$, respectively. Let~$\vtF$ denote the complex DFT (of length~$M$) of~$\vtf$; this notation applies throughout the section. For convenience, we also define $\tF_M\deq \tF_0$. After computing~$\vtF$, one can use symmetries to recover~$\vF$ (the DFT of~$\vf$); indeed, the \emph{unpacking} formula in Fourier space is
\begin{equation}\label{eqn:fourierUnpack}
F_k = \frac{1}{2}\left(\tF_k + \conj{\tF_{{\Nh}-k}}\right) - \frac{i}{2}\left(\tF_k - \conj{\tF_{{\Nh}-k}}\right)\zeta_{2\Nh}^{k}, \quad k \in\{ 0, \ldots, {\Nh}\},
\end{equation}
with inverse
\begin{equation}\label{eqn:fourierPack}
\tF_k = \frac{1}{2}\left(F_k + \conj{F_{{\Nh}-k}}\right) + \frac{i}{2}\left(F_k - \conj{F_{{\Nh}-k}}\right)\zeta_{2\Nh}^{-k}, \quad k \in \{0, \ldots, {\Nh}-1\}.
\end{equation}
Since $\vF$ is Hermitian symmetric, values of $F_k$ for $k\in \{\Nh+1,\ldots, 2\Nh-1\}$ need not be computed. For details of the derivation of \cref{eqn:fourierUnpack,eqn:fourierPack}, see \cite{PressEtAl2007}.

Complex packing enables the computation of real DFTs via complex DFTs of half the length. For hybrid dealiasing this would permit use of the complex convolution algorithms of \cite[\S2]{MuraskoBowman2024}. However, despite the practical efficiency of complex FFTs over real-to-complex/complex-to-real variants, explicit packing and unpacking in Fourier space is costly: packing-based DFT algorithms for real inputs require ``more additions than a specialized algorithm for real input data'' \cite{SorensenEtAl1987}. However, since our objective is the computation of a convolution rather than a DFT, we can perform the packing implicitly via a custom multiplication operator.

Let~$\vg\deq \seq{g}{j}{0}{2\Nh-1}$ be another real array, $\vtg$ its complex packing, and~$\vtG$ the DFT of~$\vtg$. Using \cref{eqn:fourierUnpack} and \cref{eqn:fourierPack}, define the \emph{Fourier unpacking operator} $U:\vtF\mapsto \vF$ and its inverse $U^{-1}:\vF\mapsto \vtF$. We then define the real convolution multiplication (RCM) operator $\boxdot$ to be the binary operation on~$\C^{\Nh}$ given by
\begin{equation}\label{eqn:RCM-U}
\vtF\boxdot \vtG \deq U^{-1}(U(\vtF)\odot U(\vtG)),
\end{equation}
where $\odot$ denotes element-wise multiplication. By the convolution theorem, $\vtF\boxdot \vtG$ is the DFT of the complex packing of~$\vf*\vg$. Crucially, $\vtF\boxdot \vtG$ can be computed without explicitly computing $\vF$ and $\vG$. Indeed, expanding and simplifying \cref{eqn:RCM-U} gives
\begin{equation}\label{eqn:RCM}
(\vtF \boxdot \vtG)_k \deq \tF_k\tG_k - \frac{1}{4}\left(\tF_k - \conj{\tF_{\Nh-k}}\right)\left(\tG_k - \conj{\tG_{\Nh-k}}\right)\left(1+\zeta_{{\Nh}}^{k}\right),
\end{equation}
for~$k\in\{0,\ldots, {\Nh}-1\}$. The derivation of \cref{eqn:RCM} is straightforward but lengthy; for completeness we include it in \cref{sec:RCM_derivation}.

The principal difficulty with \cref{eqn:RCM} is that it is not an element-wise operator: computing $(\vtF \boxdot \vtG)_k$ requires $\tF_k$, $\tG_k$, $\tF_{\Nh-k}$, and~$\tG_{\Nh-k}$. This imposes the constraint that at least two residues must be computed simultaneously. Using the framework of \cref{sec:recap}, our goal is to compute $(\vtF\boxdot\vtG)_{q\ell+r}$. Since $\Nh-k=qm-q\ell-r=q(m-\ell-1)+q-r$, define
$$\AA_{\ell,r}\deq \frac{1}{4}\left(\tF_{q\ell+r} - \conj{\tF_{q(m-\ell-1)+q-r}}\right)\left(\tG_{q\ell+r} - \conj{\tG_{q(m-\ell-1)+q-r}}\right)\left(1+\zeta_{qm}^{q\ell + r}\right),$$
so that
\begin{equation}\label{eqn:Hr}
(\vtF\boxdot\vtG)_{q\ell+r} = \tF_{q\ell+r}\tG_{q\ell+r}-\AA_{\ell,r}.
\end{equation}
Computing $\AA_{\ell,r}$ requires residues $r$ and~$q-r$. The cases
$r=0$ and $r=q/2$ are treated together, as they are
self-conjugate.\footnote{That is, $(q-r)\bmod q = r$ when $r=0$ or
$r=q/2$.} For $r\not\in\{0, q/2\}$, one easily verifies that
$\conj{\AA_{\ell,r}}=\AA_{m-\ell-1,q-r}$, so that on mapping $\ell\mapsto m-\ell-1$
and $r\mapsto q-r$, \cref{eqn:Hr} becomes
\begin{align}\label{eqn:Hr-conj}
(\vtF\boxdot\vtG)_{q(m-\ell-1)+q-r}= \tF_{q(m-\ell-1)+q-r}\tG_{q(m-\ell-1)+q-r}-\conj{\AA_{\ell,r}}.
\end{align}
\Cref{eqn:Hr,eqn:Hr-conj} allow simultaneous computation of residues $r$ and~$q-r$. This formulation also enables the multiplication to be performed in-place.

\subsection{Implicitly packed convolutions of more than two arrays}\label{sub:implicitly_packed_convolutions_of_more_than_two_arrays}

Although binary convolution is the most common case, some applications require convolving more than two arrays. The implicitly-packed multiplication operator extends from \cref{eqn:RCM}. Let $\vf^1, \ldots, \vf^n$ be real arrays of length $2\Nh$. For $j=1,\ldots,n$, let $\vtf^j$ be the complex packing of~$\vf^j$ (as in \eqref{eqn:complexPacking}) and let~$\vtF^j$ be the DFT of~$\vtf^j$. We first see that the RCM operator is associative:
\begin{align*}
(\vtF^1\boxdot \vtF^2)\boxdot \vtF^3 &= U^{-1}(U(\vtF^1)\odot U(\vtF^2))\boxdot \vtF^3 = U^{-1}(U(\vtF^1)\odot U(\vtF^2) \odot U(\vtF^3))\\
&=\vtF^1\boxdot U^{-1}(U(\vtF^2)\odot U(\vtF^3)) =\vtF^1\boxdot (\vtF^2\boxdot \vtF^3).
\end{align*}
By the convolution theorem, the DFT of the complex packing of~$\vf^1*\dots*\vf^n$ is then
$$U^{-1}(U(\vtF^1)\odot\dots\odot U(\vtF^n))=\vtF^1\boxdot \cdots \boxdot\vtF^n.$$
Efficient computation therefore reduces to evaluating $\vtF^1\boxdot \cdots \boxdot\vtF^n$.

This may be computed using a recursive formula. Define $\vR^1\deq \vtF^1$, and for $j\in \{2,\ldots, n\}$, define $\vR^{j} \deq \vtF^1\boxdot \cdots \boxdot\vtF^{j}$. One can compute $\vR^{j}$ using $\vR^{j-1}$ and $\tF^j$. For $k=0,\ldots, \Nh-1$, define
\begin{equation}\label{eqn:Akn}
\AA_k^j\deq \frac{1}{4}\left(R_k^{j-1}-\conj{R_{\Nh-k}^{j-1}} \right)\left( \tF_{k}^j-\conj{\tF_{\Nh-k}^j} \right)\left( 1+\zeta_{\Nh}^{k} \right).
\end{equation}
It is easy to check that $\AA_{\Nh-k}^j=\conj{\AA_{k}^j}$. This allows simultaneous computation of
\begin{equation}\label{eqn:Rnk}
R^j_k=(\vR^{j-1}\boxdot\vtF^j)_k=R_k^{j-1}\tF_k^{j}-\AA_k^j,
\end{equation}
and
\begin{equation}\label{eqn:RnMk}
R^j_{\Nh-k}=R_{\Nh-k}^{j-1}\tF_{\Nh-k}^{j}-\conj{\AA_{\Nh-k}^j}.
\end{equation}

\subsection{Comparison with explicit packing}\label{sub:comparison_of_implicit_packing_with_explicit_packing}

We compare the floating-point operation (flop) counts of implicit and explicit packing (using \cref{eqn:fourierUnpack} and \cref{eqn:fourierPack} directly). Consider convolving $n$ arrays $\vtF^{1}, \ldots, \vtF^n$, each of length $M$. For simplicity, we assume that $M$ is even.

For explicit packing, let $k=1,\ldots, M/2-1$ (the indices $k=0, M/2$ contribute only $O(n)$ terms). For each $j=1,\ldots, n$, compute
$$A_k^j\deq \frac{1}{2}\left(\tF_{k}^j+\conj{\tF_{M-k}^j}\right),\quad B_k^j\deq \frac{i}{2}\left(\tF_{k}^j-\conj{\tF_{M-k}^j}\right)\zeta_{2M}^k.$$
Computing $A_k^j$ requires 1 complex addition (2 flops) and 1 real-complex multiplication (2 flops). With the constant factor $\frac{i}{2}\zeta_{2M}^k$ precomputed, $B_k^j$ requires 1 complex addition (2 flops) and 1 complex multiplication (6 flops). From \cref{eqn:fourierUnpack},
$$F_k^j=A_k^j-B_k^j,\quad F_{M-k}^j=\conj{A_k^j+B_k^j}.$$
Assuming conjugation is free, this adds two complex additions (4 flops). Computing $F_k^j$ and $F_{M-k}^j$ thus requires 16 flops per pair. With approximately $M/2$ pairs and $n$ arrays, computing $\vF^{1}, \ldots, \vF^n$ requires approximately $8nM$ flops. Multiplying these arrays requires $(n-1)M$ complex multiplications ($6(n-1)M$ flops). Applying \cref{eqn:fourierPack} to the product requires another $8M$ flops. The total is approximately $2M(7n+1)$ flops.

For implicit packing, each step $j=2,\ldots, n$ computes $\AA^j_k$ (as in \cref{eqn:Akn}) using 2 complex multiplications (12 flops) and 2 complex additions (4 flops). Computing $R^j_k$ and $R^j_{M-k}$ via \cref{eqn:Rnk,eqn:RnMk} adds 2 more complex multiplications (12 flops) and 2 more complex additions (4 flops), for 32 flops per pair. With approximately $M/2$ pairs processed $n-1$ times, the total is approximately $16M(n-1)$ flops.

For a binary convolution ($n=2$), implicit packing requires approximately half as many flops as explicit packing. However, its flop count scales worse with $n$: for $n>9$, it exceeds that of explicit packing. Of course, this analysis uses the recursive formula of \cref{sub:implicitly_packed_convolutions_of_more_than_two_arrays}; it is possible that a more optimal formula may reduce the count for large $n$.

\section{Multidimensional real convolutions}\label{sec:multidim} To compute multidimensional real convolutions, we use the decomposition of
\cite{BowmanRoberts2011,RobertsBowman2018}.

An $n$-dimensional complex convolution is conventionally computed using complex
FFTs of size~$N_1\times\ldots\times N_n$. Instead, we decompose the
$n$-dimensional convolution recursively into $\prod_{i=2}^{n}N_i$ FFTs in the first dimension, followed by $N_1$ convolutions of
dimension $n-1$, and finally $\prod_{i=2}^{n}N_i$ inverse FFTs in the first
dimension. This reduces a multidimensional convolution to one-dimensional convolutions at the innermost level, permitting buffer reuse across subconvolutions and substantially reducing memory usage. For more details on this algorithm, see \cite[\S5]{MuraskoBowman2024}.

This decomposition applies directly to the algorithms of \cref{sec:symmetries}. The only difference is that the outermost DFT must handle real values; all subconvolutions are complex\footnote{This contrasts with the Hermitian-symmetric case, where only the innermost FFTs need to be complex-to-real/real-to-complex.} and use the hybrid-dealiased convolutions of \cite{MuraskoBowman2024}.

In contrast, applying the decomposition to the implicit packing algorithms of \cref{sec:rcm} is much more involved. Let $M_1,\ldots, M_n\in \N$. For this section we use indices $k_1,\dots, k_n$ with $0\leq k_j < M_j$ for $j=1,\ldots, n-1$, and write $k_j^*\deq (M_j-k_j)\bmod M_j$. The multi-index is $\vk\deq (k_1,\ldots, k_n)$ with conjugate $\vk^*\deq (k_1^*,\ldots,k_n^*)$.

Unless stated otherwise, all $n$-dimensional arrays are complex and of size~$M_1\times \dots\times  M_n$. For any such array~$\vA$ and $0\leq k_1<M_1$, define the $(n-1)$-dimensional slice~$\vA^{k_1}$ by
\begin{equation}\label{eqn:slice}
A^{k_1}_{(k_2, \ldots, k_n)} \deq A_{k_1,k_2,\ldots, k_n}.
\end{equation}

Let $\vf$ and~$\vg$ be real $n$-dimensional arrays of size
$M_1\times \dots\times  M_{n-1}\times 2M_n$. Their complex packings $\vtf$
and~$\vtg$ are obtained by packing along the last dimension:
\begin{equation}\label{eqn:multidimComplexPacking}
\tf_{j_1,\ldots, j_n}\deq f_{j_1,\ldots, 2j_n}+i f_{j_1,\ldots, 2j_n+1}, \quad \tg_{j_1,\ldots, j_n}\deq g_{j_1,\ldots, 2j_n}+i g_{j_1,\ldots, 2j_n+1},
\end{equation}
where $0\leq j_i < M_i$ for $i=1,\ldots, n$. Let $\vtF$
and~$\vtG$ denote the $n$-dimensional DFTs of~$\vtf$ and~$\vtg$, respectively. It is straightforward to show that the RCM operator \cref{eqn:RCM} generalizes to
\begin{equation}\label{eqn:RCM-ndim}
(\vtF \boxdot_n \vtG)_\vk \deq \tF_
\vk\tG_\vk - \frac{1}{4}\left(\tF_\vk - \conj{\tF_{\vk^*}}\right)\left(\tG_\vk - \conj{\tG_{\vk^*}}\right)\left(1+\zeta_{{M_n}}^{k_n}\right).
\end{equation}

\Cref{eqn:RCM-ndim} shows why implicit packing complicates the multidimensional decomposition: computing $(\vtF \boxdot_n \vtG)_\vk$ requires elements indexed by both $\vk$ and~$\vk^*$. To handle this, we define a generalized operator. For any $n$-dimensional arrays $\vA$, $\vB$, $\vC$, and~$\vD$, let $\cC_n(\vA,\vB,\vC,\vD)$ be the array with entries
$$(\cC_n(\vA,\vB,\vC,\vD))_{\vk}={A}_
{\vk}{B}_{\vk} - \frac{1}{4}\left({A}_{\vk} - \conj{{C}_{\vk^*}}\right)\left({B}_{\vk} - \conj{{D}_{\vk^*}}\right)\left(1+\zeta_{{M_n}}^{k_n}\right).$$
One easily verifies that
\begin{align*}
(\cC_n(\vA,\vB,\vC,\vD))_{(k_1,k_2,\ldots,k_n)}
&=(\cC_{n-1}(\vA^{k_1},\vB^{k_1},\vC^{k_1^*},\vD^{k_1^*}))_{(k_2,\ldots,k_n)}.
\end{align*}
By construction, $(\vtF \boxdot_n \vtG)=\cC_n(\vtF,\vtG,\vtF,\vtG)$. For $0\leq k_1 \leq M_1/2$, we compute $\vtF^{k_1}$, $\vtF^{k_1^*}$, $\vtG^{k_1}$, and
$\vtG^{k_1^*}$, which yields both
\begin{equation*}
(\vtF \boxdot_n \vtG)_{(k_1,k_2,\ldots, k_n)}=(\cC_{n-1}(\vtF^{k_1},\vtG^{k_1},\vtF^{k_1^*},\vtG^{k_1^*}))_{(k_2,\ldots, k_n)},
\end{equation*}
and
\begin{equation*}
(\vtF \boxdot_n \vtG)_{(k_1^*,k_2,\ldots, k_n)}=(\cC_{n-1}(\vtF^{k_1^*},\vtG^{k_1^*},\vtF^{k_1},\vtG^{k_1}))_{(k_2,\ldots, k_n)}.
\end{equation*}
These contributions are handled together. Applying this recursively reduces $n$-dimensional implicitly packed convolutions to a sequence of 1-dimensional implicitly packed convolutions.

\begin{figure}[htbp]
\begin{center}
\includegraphics{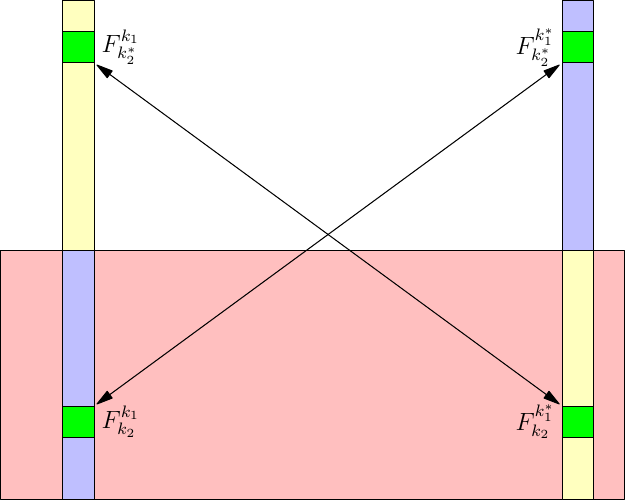}
\end{center}
\caption{Relative positions of conjugate modes in 2D implicit packing. The DFT is first computed along the horizontal dimension. For each $k_1$, two 1D implicitly packed convolutions are computed along the vertical dimension: the blue array (bottom left and top right) serves as input to the first convolution and the yellow array (top left and bottom right) to the second.}
\label{fig:2Drcm}
\end{figure}

Intuitively, the first two arguments of~$\cC_{n-1}$ correspond to the \emph{bottom half} of the array and the last two to the \emph{top half}. \Cref{fig:2Drcm} illustrates this in two dimensions.

Pseudocode for one- and two-dimensional implicitly packed convolutions is given in \cref{alg:convolveRCMone,alg:convolveRCMtwo}, both using the common multiplication routine \cref{alg:multBinaryRCM}. For simplicity, the pseudocode considers only explicit zero padding, but may be extended to the hybrid dealiasing algorithms of \cite{MuraskoBowman2024}.

\begin{multalg}[htbp]
\hfill
\begin{minipage}[t]{0.49\textwidth}
\centering
\begin{algorithm}[H]
\algsize
\caption{\convolveRCMone{} returns an explicitly dealiased and implicitly packed one-dimensional convolution of two real arrays~$\vf$ and~$\vg$. It makes use of routines \zeropad, which performs the necessary explicit zero padding, as well as \complexpack{} and \complexunpack{}, which pack and unpack the real data as in \cref{eqn:complexPacking}. The multiplication is done using \cref{alg:multBinaryRCM}.}
\begin{algorithmic}\label{alg:convolveRCMone}
\setlength{\baselineskip}{\algspace\baselineskip}
\Input{$\vf$, $\vg$}
\State{$\vtf\gets \complexpack(\zeropad(\vf))$}
\State{$\vtg\gets \complexpack(\zeropad(\vg))$}
\State{$\vtF \gets \fft{\vtf}$}
\State{$\vtG \gets \fft{\vtg}$}
\State{$\multBinaryRCM(\vtF,\vtG,\vtF,\vtG,\True)$}
\State{$\vth \gets \ifft{\vtF}$}
\State{$\vh \gets \complexunpack(\vth)$}\\
\Return{$\vh$}
\end{algorithmic}
\end{algorithm}

\vspace{-1.5em}

\begin{algorithm}[H]
\algsize
\caption{\multBinaryRCM\ is an in-place implicit packing multiplication routine. It takes four length~$M$ one-dimensional complex arrays~$\vA,\vB,\vC,\vD$ (where we assume~$M$ is even), and a boolean value \first. Output is written to~$\vA$ and~$\vC$. This routine makes use of pre-computed values~$Z_k=(1+\zeta_M^k)/4$ where~$k=1,\ldots, M/2-1$. }
\begin{algorithmic}\label{alg:multBinaryRCM}
\setlength{\baselineskip}{\algspace\baselineskip}
\Input{$\vA$,$\vB$,$\vC$,$\vD$, \first}
\If{\first}
  \State{$\alpha \gets \frac{1}{2} \left(A_0 - \conj{C_0}\right)\left(B_0 - \conj{D_0}\right)$}
  \State{$A_0\gets A_0 B_0 - \alpha$}
  \If{\Not($\vA = \vC$ \Or{}~$\vB= \vD$)}
    \State{$C_0\gets C_0 D_0 - \conj{\alpha}$}
  \EndIf
\EndIf
\For{$k=1,\ldots,M/2-1$}
  \State{$\alpha\gets \left( A_k - \conj{C_{M-k}} \right)\left( B_k - \conj{D_{M-k}} \right) Z_k$}
  \State{$A_k\gets A_k B_k - \alpha$}
  \State{$C_{M-k}\gets C_{M-k} D_{M-k} - \conj{\alpha}$}
\EndFor
\State{$A_{M/2}\gets A_{M/2} B_{M/2}$}
\end{algorithmic}
\end{algorithm}

\end{minipage}
\hfill
\begin{minipage}[t]{0.49\textwidth}
\centering

\begin{algorithm}[H]
\algsize
\caption{\convolveRCMtwo{} returns an explicitly dealiased and implicitly packed two-dimensional convolution of two real arrays~$\vf$ and~$\vg$. It makes use of routines \zeropad, which explicitly zero pads the arrays to size~$M_1\times 2M_2$, as well as \complexpack{} and \complexunpack{}, which pack and unpack the real data along the last dimension as in \cref{eqn:multidimComplexPacking}. The routine \fftx{} computes FFTs along the first dimension. We use the notation defined in \cref{eqn:slice} for one-dimensional slices. The multiplication is done using \cref{alg:multBinaryRCM}.}
\begin{algorithmic}\label{alg:convolveRCMtwo}
\setlength{\baselineskip}{\algspace\baselineskip}
\Input{$\vf$, $\vg$}
\State{$\vtf\gets \complexpack(\zeropad(\vf))$}
\State{$\vtg\gets \complexpack(\zeropad(\vg))$}
\State{$\vtF \gets \fftx(\vtf)$}
\State{$\vtG \gets \fftx(\vtg)$}
\For{$k_1=0,\ldots,M_1/2-1$}
  \State{$k_1^*\gets M_1-k_1$}
  \If{$k_1=0$}
    \State{$k_1^*\gets M_1/2$}
  \EndIf
  \State{$\vtF^{k_1}\gets \ffty(\vtF^{k_1})$}
  \State{$\vtG^{k_1}\gets \ffty(\vtG^{k_1})$}
  \State{$\vtF^{k_1^*}\gets \ffty(\vtF^{k_1^*})$}
  \State{$\vtG^{k_1^*}\gets \ffty(\vtG^{k_1^*})$}
  \If{$k_1=0$}
    \State{$\multBinaryRCM(\vtF^{k_1},\vtG^{k_1},\vtF^{k_1},\vtG^{k_1}, \True)$}
    \State{$\multBinaryRCM(\vtF^{k_1^*},\vtG^{k_1^*},\vtF^{k_1^*},\vtG^{k_1^*}, \True)$}
  \Else
    \State{$\multBinaryRCM(\vtF^{k_1},\vtG^{k_1},\vtF^{k_1^*},\vtG^{k_1^*}, \True)$}
    \State{$\multBinaryRCM(\vtF^{k_1^*},\vtG^{k_1^*},\vtF^{k_1},\vtG^{k_1}, \False)$}
  \EndIf
  \State{$\vtF^{k_1}\gets \iffty(\vtF^{k_1})$}
  \State{$\vtF^{k_1^*}\gets \iffty(\vtF^{k_1^*})$}
\EndFor
\State{$\vth \gets \ifftx(\vtF)$}
\State{$\vh \gets \complexunpack(\vth)$}\\
\Return{$\vh$}
\end{algorithmic}
\end{algorithm}
\end{minipage}
\end{multalg}

\section{Numerical results}\label{sec:results}

To benchmark the above algorithms, we follow \cite[\S7]{MuraskoBowman2024}, tuning the parameters empirically by scanning over the FFT size~$m$, the number~$D$ of residues computed simultaneously, and the choice between in-place and out-of-place FFTs. For multidimensional convolutions we assume that the parameters for each dimension can be optimized independently. We consider only power-of-two sizes (optimal for FFTs), even though hybrid dealiasing performs well for arbitrary sizes (see \cite[Fig.~17]{MuraskoBowman2024}).

The benchmarks use \fftwpp~\cite{fftwpp}, which contains vectorized and parallelized C++ implementations of these algorithms. We benchmark on an Intel i9-12900K processor (5.2 GHz, 8 performance cores) with an ASUS ROG Strix Z690-F motherboard and 128 GB DDR5 memory (5 GHz), compiled with GCC 15.2.1 using flags {\tt -Ofast -fomit-frame-pointer -fstrict-aliasing -ffast-math}. The underlying FFTs use FFTW 3.3.10 \cite{fftw,FrigoJohnson2005} under Fedora 43, with multithreading via OpenMP.

We measure the time to compute an unnormalized\footnote{By unnormalized convolution, we mean that the inverse FFTs are not normalized.} in-place dealiased binary convolution using the {\tt chrono::nanoseconds} clock from the C++ standard library. Each size is timed for at least 5 seconds and until at least 20 samples are collected. We plot the median times from each of these tests. We compare four algorithms:
\begin{itemize}
\item \explicit uses real-to-complex/complex-to-real FFTs and explicit padding (the standard approach). In one dimension, both in-place and out-of-place FFTs are benchmarked and the minimum time is plotted. In two and three dimensions, in-place multidimensional FFTs are used for all sizes.
\item \hybridCS uses the conjugate-symmetry algorithm of \cref{sec:symmetries}. In two and three dimensions, it employs the multidimensional decomposition of \cref{sec:multidim} with complex FFTs from \cite{MuraskoBowman2024} for the inner transforms.
\item \hybridRCM uses the implicit packing algorithm of \cref{sec:rcm} with complex hybrid dealiasing from \cite{MuraskoBowman2024}. In two and three dimensions, it employs the multidimensional decomposition of \cref{sec:multidim}.
\item \explicitRCM is a simplified variant of \hybridRCM: it uses complex packing and the RCM operator but only explicit zero padding for dealiasing.
\end{itemize}

\subsection{Convolution benchmarks in one dimension}\label{sub:one_dimensional_convolutions}
\Cref{fig:timingsr1-T1} plots median execution times (normalized to $L\log_2 L$) for one-dimensional in-place convolutions of~$L$ real words ($M=2L$) on a single thread. Both \hybridCS and \explicitRCM are, on average, slightly faster than \explicit. On the other hand, the \hybridRCM algorithm outperforms \explicit at all sizes and is significantly faster for large convolutions, with a maximum speedup factor of approximately 1.4.

\Cref{fig:timingsr1-T8} plots normalized times for the same convolutions parallelized over eight threads. For large sizes, \explicit is outperformed by all other algorithms. Again, the \hybridRCM algorithm is the clear winner, with a maximum speedup factor of approximately 3.9 over \explicit. Notably, \explicitRCM performs nearly as well as \hybridRCM for many large sizes, although its performance degrades for the largest convolutions.

\begin{figure}[tbhp]
\begin{minipage}{0.49\linewidth}
\Figure{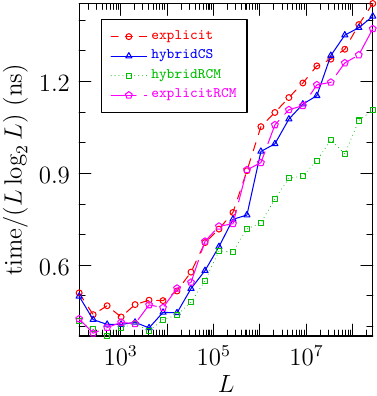}{\benchCaptionOne{1 thread}}
\end{minipage}
\,
\begin{minipage}{0.49\linewidth}
\Figure{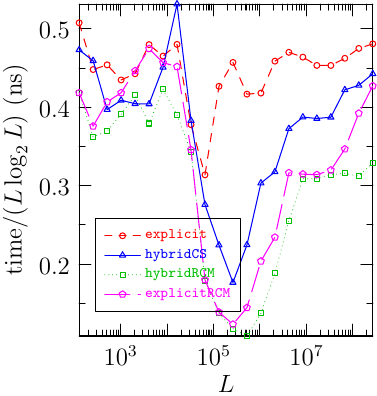}{\benchCaptionOne{8 threads}}
\end{minipage}
\end{figure}

Although \hybridRCM and \explicitRCM often achieve similar timings, we note that they employ very different algorithms: for these large cases, the optimal parameters for \hybridRCM always use the complex inner-loop optimization of \cite[\S 6.1]{MuraskoBowman2024}. Under multithreading, most of the speedup from implicit packing appears to stem from the ability to use complex FFTs rather than real-to-complex/complex-to-real variants.

\subsection{Convolution benchmarks in two and three dimensions}\label{sub:convolution_benchmarks_in_two_and_three_dimensions}
\Cref{fig:timingsr2-T1,fig:timingsr2-T8} plot normalized times for two-dimensional real convolutions of size $L\times L$ padded to $2L\times 2L$, on one thread and eight threads respectively. As in \cite{MuraskoBowman2024}, \explicit is outperformed by all other algorithms. On both one and eight threads, \hybridRCM matches or exceeds all others. On one thread, \hybridCS now outperforms \explicitRCM at all sizes; on eight threads, \explicitRCM usually outperforms \hybridCS. We also observe degraded performance of \explicitRCM at the largest size, consistent with the one-dimensional eight-thread results in \cref{fig:timingsr1-T8}.

\begin{figure}[tbhp]
\begin{minipage}{0.49\linewidth}
\Figure{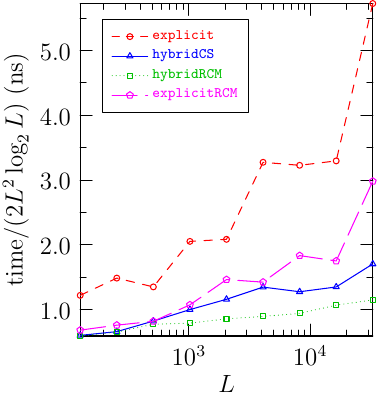}{\benchCaptionTwo{1 thread}}
\end{minipage}
\,
\begin{minipage}{0.49\linewidth}
\Figure{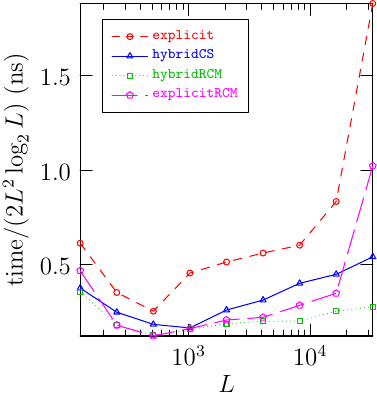}{\benchCaptionTwo{8 threads}}
\end{minipage}
\end{figure}

\Cref{fig:timingsr3-T1,fig:timingsr3-T8} plot normalized times for three-dimensional real convolutions of size $L\times L\times L$ padded to $2L\times 2L\times 2L$, on one thread and eight threads respectively. Again, \explicit is outperformed by all other algorithms, and \hybridRCM is the clear winner on both one and eight threads (except at the smallest benchmark). Both \hybridCS and \explicitRCM are comparable here.

\begin{figure}[tbhp]
\begin{minipage}{0.49\linewidth}
\Figure{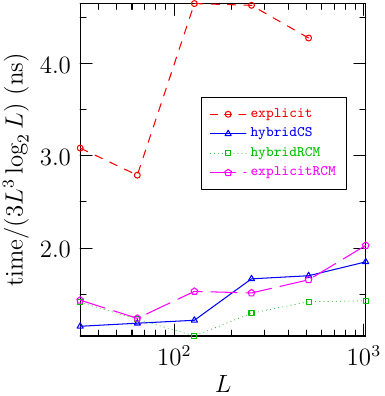}{\benchCaptionThree{1 thread}}
\end{minipage}
\,
\begin{minipage}{0.49\linewidth}
\Figure{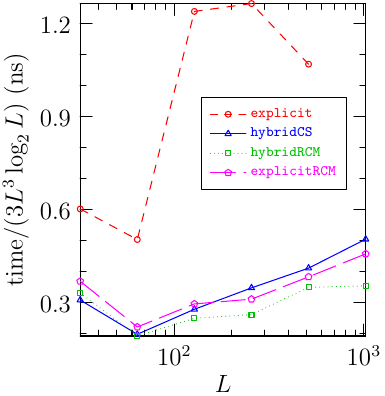}{\benchCaptionThree{8 threads}}
\end{minipage}
\end{figure}

\subsection{Discussion}\label{sub:discussion}
The results of \cref{sec:results} indicate that \hybridRCM consistently
outperforms all other tested algorithms. Nevertheless, both \explicitRCM and \hybridCS have distinct merits.

Benchmarking \explicitRCM is motivated by the difficulty of implementing efficient hybrid-dealiased convolutions. While such routines are available in the high-performance C++ library \fftwpp \cite{fftwpp}, some applications cannot integrate this library (for instance, no GPU implementation of hybrid dealiasing currently exists). Implementing implicitly packed convolutions with explicit dealiasing and multidimensional decomposition is considerably simpler. Our benchmarks suggest this is worthwhile: in several cases \explicitRCM performs nearly as well as \hybridRCM, particularly under multithreading.

The \hybridCS algorithm also offers advantages over \hybridRCM. First, its multiplication operator generalizes easily. As discussed in
\cref{sub:implicitly_packed_convolutions_of_more_than_two_arrays}, the RCM operator \cref{eqn:RCM} can be extended to convolutions of $n\geq 2$ real arrays, but deriving and implementing a multiplication operator for implicit packing may be difficult for more complicated generalized convolutions. Moreover, even when derivable, an implicitly packed multiplier may not be efficient (as shown in \cref{sub:comparison_of_implicit_packing_with_explicit_packing} when convolving more than 9 arrays at once). Second, some applications require the Fourier transforms of the inputs alongside the convolution; these are computed explicitly by \hybridCS but not by \hybridRCM.

An alternative approach that combines \hybridCS and \hybridRCM may be worth investigating: computing residue contributions is more involved in
\hybridRCM than in \hybridCS. As discussed in \cref{sec:rcm}, for residues $r\not\in \{0,q/2\}$ we must compute residues $r$ and $q-r$ simultaneously, increasing the required work memory. Since this is particularly significant in higher dimensions, one could follow the \hybridCS algorithm of \cref{sec:symmetries} but use implicit packing for the $0$ and~$q/2$ residue contributions. This would allow each residue to be computed independently.

\section{Conclusion and future work}\label{sec:conclusion}
The original motivation for implicit dealiasing (the predecessor to hybrid dealiasing) was improving the efficiency of pseudospectral simulations of partial differential equations, which rely on convolutions of Hermit\-ian-symmetric arrays. Extending hybrid dealiasing to real-valued inputs substantially broadens the range of applicable problems. In particular, these algorithms may accelerate convolutional neural networks (CNNs).

In a CNN, input data is convolved with a small kernel of learned parameters. These convolutions are typically computed directly or via the Winograd method \cite{LavinGray2016}, though FFT-based approaches have also been explored \cite{MathieuEtAl2014,HighlanderRodriguez2016,LinYao2019,ChitsazEtAl2020}. FFT-based convolutions become competitive for large kernels, large inputs, and high-dimensional problems. For small kernels, aliasing errors affect only a narrow boundary region, and it is common to use ``valid'' padding, discarding outputs near the boundary after convolution. The primary purpose of zero padding is then to match the kernel size to the input arrays. Since kernels are typically much smaller than inputs, this entails substantial zero padding.

One direction for future work is modifying the hybrid dealiasing algorithms to handle arrays of different sizes, enabling implicit padding of the kernel. Following the derivation of \cref{sec:recap}, one could allow the arrays to use different values of $p$. Specifically, to convolve an array $\vf$ of length $L_f$ with a kernel $\vg$ of length $L_g\leq L_f$, both zero padded to length $M\geq L_f$: first find the optimal FFT size $m$ for $\vf$, then explicitly pad $\vf$ and $\vg$ to lengths $p_fm$ and $p_g m$ respectively, where $p_f\deq\ceil{L_f/m}$ and $p_g\deq \ceil{L_g/m}$. Both arrays are then implicitly padded and transformed to length $qm$ with $q\deq \ceil{M/m}$. Since both arrays share the same $q$ and $m$, the convolution can still be computed one residue at a time. Given that the required zero padding scales exponentially with dimension, this approach may benefit high-dimensional convolutions.

However, implicit zero padding of the kernel may be unnecessary in some applications. In a CNN, zero padding and transforming the kernel is needed only when the kernel is updated. This does not occur during inference, and only occurs in training after each mini-batch. The padded, transformed kernel can therefore be precomputed and incorporated into the multiplication routine (though this must be done for each input size). If input size variability is low, the padding and transform cost will likely be amortized.

Future work could also extend hybrid dealiasing to block-based convolution algorithms such as overlap-add or overlap-save. These methods partition the input into overlapping blocks, convolve each with the kernel, and combine the results. They reduce FFT length (and hence padding/transform costs), particularly for very large inputs. The block size introduces an additional optimization parameter; it would be interesting to study how different block sizes affect the optimal hybrid dealiasing parameters.

Finally, a GPU implementation of hybrid dealiasing is a priority. To date, these algorithms have been implemented and benchmarked only on CPUs, whereas standard CNN training and inference pipelines are typically GPU-accelerated.

As noted in \cref{sub:discussion}, a plausible starting point is implementing \explicitRCM on the GPU, using implicit packing with explicit dealiasing. This would leverage existing GPU FFT-based convolution techniques. The RCM operator, though more complex than Hadamard multiplication, could be implemented as a custom kernel. The results of \cref{sec:results} suggest that even without hybrid dealiasing, a GPU version of \explicitRCM may be enough to outperform standard convolution methods.

\appendix

\section{Inner loop for real arrays}\label{sec:realInner}
The algorithms presented so far support any FFT size $m\in \N$, but efficiency degrades when $m$ is much smaller than the input arrays. As $m$ decreases, both $p$ and $q$ increase, making the preprocessing and postprocessing in \cref{eqn:forwardr,eqn:backwardr} considerably more expensive. The inner loop optimization performs these steps using FFTs of size $p$. This was described for the complex case in \cite[\S6.1]{MuraskoBowman2024} (and is therefore used in the implicit packing algorithm of \cref{sec:rcm}). Here we develop it for the algorithm of \cref{sec:symmetries}. Pseudocode appears in \cite[Algorithms 22 and 23]{Murasko2024_thesis}.

We begin by redefining $q$: let $n\deq \ceil{M/(pm)}$ and set $q\deq np$. We reindex
$$r=un+v, \quad u=0,\ldots, p-1,\quad v=0,\ldots, n-1.$$
From \cref{eqn:forwardr}, the forward transform is
\begin{equation*}
  F_{q\ell+un+v} =\sum_{s=0}^{m-1}\ze{m}{\ell s}\ze{qm}{(un+v)s}
  \sum_{t=0}^{p-1}\ze{p}{ut}\left(\ze{q}{vt}f_{tm+s}\right),
\end{equation*}
and from \cref{eqn:backwardr}, the inverse is
\begin{equation*}
  f_{tm+s} = \frac{1}{qm}\sum_{v=0}^{n-1}\ze{q}{-tv}
  \sum_{u=0}^{p-1}\ze{p}{-tu}\ze{qm}{-s(un+v)}
  \sum_{\ell=0}^{m-1}\ze{m}{-s\ell}F_{q\ell+un+v}.
\end{equation*}
Following \cref{sub:oneRes}, we consider one residue at a time (with $v$ replacing $r$). Define $\vh_v\deq\{h_{v,j}\}_{j=0}^{pm-1}$ where
\begin{equation}\label{eqn:hv}
h_{v,tm+s}\deq\zeta_{q}^{-tv}\sum_{u=0}^{p-1}\zeta_{p}^{-tu}\zeta_{qm}^{-s(un+v)}\sum_{\ell=0}^{m-1}\zeta_{m}^{-s\ell}F_{q\ell+un+v},
\end{equation}
so that $\vf=\frac{1}{qm}\sum_{v=0}^{n-1}\vh_{v}$. As in \cref{eqn:hrHerm}, the Hermitian symmetry $\conj{\vh_{v}}=\vh_{n-v}$ gives
\begin{equation}\label{eqn:real-accumulationInner}
\vf=\frac{1}{qm}\left(\vh_{0}+2\sum_{v=1}^{\ceil{n/2}-1}\Re \vh_{v}+\vh_{n/2}\right),
\end{equation}
where $\vh_{n/2} \equiv \vector{0}$ if $n$ is odd. As in \cref{sub:conjugate}, there are three cases.

\subsection{Case \texorpdfstring{$\vh_0$}{v=0}}
This case reduces essentially to the algorithm of \cref{sub:conjugate}. When $v=0$, the forward transform becomes
\begin{equation}\label{eqn:Fv0}
F_{q\ell+un} =\sum_{s=0}^{m-1}\zeta_{m}^{\ell s} \zeta_{qm}^{uns}\sum_{t=0}^{p-1}\zeta_{p}^{ut}f_{tm+s}.
\end{equation}
The sum in $t$ is computed with a real-to-complex FFT of length $p$. For the sum in $s$, note that \cref{eqn:Fv0} has the same form as \cref{eqn:forwardr} with $r=un$. We therefore apply the algorithm from \cref{sub:conjugate}: when $u=0$ we use a real-to-complex FFT of length~$m$.\footnote{For simplicity, our implementation in \fftwpp \cite{fftwpp} uses a complex FFT for this part. The inner loop is typically useful when $p$ is large, so this does not significantly affect performance.} When
$u=1,\ldots, \ceil{p/2}-1$, we use complex FFTs of length~$m$. Finally, when $u=p/2$ we use a complex FFT of length~$m/2$ (requiring $m$ even).

\subsection{Case \texorpdfstring{$\vh_v$, $v\in\{1,\ldots,\ceil{n/2}-1\}$}{v=1,...,ceil(n/2)-1}}

From \cref{eqn:real-accumulationInner}, the contributions of $\vh_v$ and~$\vh_{n-v}$ are obtained simultaneously. We therefore use complex FFTs without sacrificing efficiency, following the complex inner loop algorithm of \cite[\S 6.1]{MuraskoBowman2024}: the forward transform uses $m$ FFTs of length $p$ followed by $p$ complex FFTs of length~$m$, and the inverse uses $p$ complex FFTs of length~$m$ followed by $m$ complex FFTs of length~$p$.

\subsection{Case \texorpdfstring{$\vh_{n/2}$}{v=n/2}}

This case arises only when $n$ is even. The forward transform becomes
\begin{equation}\label{eqn:Fvn2}
F_{q\ell+un+n/2} =\sum_{s=0}^{m-1}\zeta_{m}^{\ell s} \zeta_{qm}^{(un+n/2)s}\sum_{t=0}^{p-1}\zeta_{p}^{ut}\zeta_{2p}^{t}f_{tm+s},
\end{equation}
and \cref{eqn:hv} becomes
\begin{equation*}
h_{n/2,tm+s}\deq\zeta_{2p}^{-t}\sum_{u=0}^{p-1}\zeta_{p}^{-tu}\zeta_{qm}^{-s(un+n/2)}\sum_{\ell=0}^{m-1}\zeta_{m}^{-s\ell}F_{q\ell+un+n/2}.
\end{equation*}
Consider the DFT of length $p$ in \cref{eqn:Fvn2}:
\begin{equation}\label{eqn:pDFTn2}
\sum_{t=0}^{p-1}\zeta_{p}^{ut}\zeta_{2p}^{t}f_{tm+s}.
\end{equation}
Setting $p=1$ in \cref{eqn:resq2}, the sum over $s$ has the same form as the sum over $t$ in \cref{eqn:pDFTn2}. If $p$ is even, we can therefore use the algorithm for $r=q/2$ from \cref{ssub:q2} to compute \cref{eqn:pDFTn2} using complex FFTs of size $p/2$.

Let $\P\deq p/2$, and reindex $t$ and $u$ as follows:
$$t=\alpha \P +\beta, \ \ \alpha\in\{0,1\}, \ \beta\in\{0, \ldots,\P-1\}.$$
We need only compute \cref{eqn:Fvn2} for even values of $u$.
For $\gamma=0,\ldots, \P-1$ we obtain
\begin{equation}\label{eqn:Fn2gamma}
F_{q\ell+2\gamma n+n/2}=\sum_{s=0}^{m-1}\zeta_{m}^{\ell s}\zeta_{qm}^{(2\gamma n+n/2)s}\sum_{\beta=0}^{\P-1} \zeta_{\P}^{\gamma\beta}\zeta_{2p}^{\beta}\left[f_{\beta m+s}+if_{(\P+\beta)m+s}\right].
\end{equation}
Next we handle the inverse transform. As in \cref{eqn:wb}, define
\begin{equation*}
w_{s,\beta}\deq \zeta_{2p}^{-\beta}\sum_{\gamma=0}^{\P-1}\zeta_{\P}^{-\beta\gamma}\zeta_{qm}^{-s(2\gamma n+n/2)}\sum_{\ell=0}^{m-1}\zeta_{m}^{-s\ell}F_{q\ell+2\gamma n+n/2}.
\end{equation*}
This array is computed using $\P$ complex FFTs of length $m$, followed by $m$ complex FFTs of length~$\P$. We then have
\begin{equation*}
h_{n/2,\beta m +s}=2\Re\left(w_{s,\beta}\right), \quad h_{n/2,(\P+\beta)m+s}=2\Im\left(w_{s,\beta}\right),
\end{equation*}
by an argument analogous to that underlying \cref{eqn:accumulateq2}.

Thus the $v=n/2$ residue requires only $p/2$ DFTs of length~$m$ (instead of $p$), so complex FFTs of length~$m$ in \cref{eqn:Fn2gamma} incur no loss in efficiency.

\section{Derivation of the RCM operator}\label{sec:RCM_derivation}
Let $\vf$ and~$\vg$ be real arrays of length $N=2\Nh$, and let $\vtf$ and~$\vtg$ be their complex packings (\cref{eqn:complexPacking}). We derive \cref{eqn:RCM} from \cref{eqn:RCM-U}, identifying $\tF_\Nh\deq \tF_0$ and $\tG_\Nh\deq \tG_0$. Unpacking in Fourier space via \cref{eqn:fourierUnpack},
\begin{align*}
&(U(\vtF)\odot U(\vtG))_k\\
&= \frac{1}{4}\left[\left(\tF_k + \conj{\tF_{\Nh-k}}\right) - i\left(\tF_k - \conj{\tF_{\Nh-k}}\right)\zeta_{2\Nh}^{k}\right]\left[\left(\tG_k + \conj{\tG_{\Nh-k}}\right) - i\left(\tG_k - \conj{\tG_{\Nh-k}}\right)\zeta_{2\Nh}^{k}\right]\\
&=\frac{1}{4}\left[\tF_k\tG_k + \conj{\tF_{\Nh-k}}\tG_k + \tF_k\conj{\tG_{\Nh-k}}+\conj{\tF_{\Nh-k}}\conj{\tG_{\Nh-k}}- 2i\left(\tF_k\tG_k -\conj{\tF_{\Nh-k}}\conj{\tG_{\Nh-k}}\right)\zeta_{2\Nh}^{k}-\right.\\
&\qquad \left. \left(\tF_k\tG_k -\conj{\tF_{\Nh-k}}\tG_k -\tF_k\conj{\tG_{\Nh-k}}+\conj{\tF_{\Nh-k}}\conj{\tG_{\Nh-k}}\right)\zeta_{\Nh}^{k}\right],
\end{align*}
for all $k = 0, \ldots, \Nh$. Define
\begin{align*}
A_k &\deq \tF_k\tG_k + \conj{\tF_{\Nh-k}}\tG_k + \tF_k\conj{\tG_{\Nh-k}}+\conj{\tF_{\Nh-k}}\conj{\tG_{\Nh-k}},\\
B_k&\deq \tF_k\tG_k -\conj{\tF_{\Nh-k}}\conj{\tG_{\Nh-k}},\\
C_k&\deq \tF_k\tG_k -\conj{\tF_{\Nh-k}}\tG_k -\tF_k\conj{\tG_{\Nh-k}}+\conj{\tF_{\Nh-k}}\conj{\tG_{\Nh-k}},
\end{align*}
so that
\begin{equation}\label{eqn:modmultABC}
(U(\vtF)\odot U(\vtG))_k =\frac{1}{4}\left(A_k- 2iB_k\zeta_{2\Nh}^{k} - C_k\zeta_{\Nh}^{k}\right).
\end{equation}
Note the symmetries
$$A_k = \conj{A_{\Nh-k}},\quad B_k = -\conj{B_{\Nh-k}},\quad C_k = \conj{C_{\Nh-k}}.$$
Using these, we compute
\begin{equation}\label{eqn:modmultABCconj}
\conj{(U(\vtF)\odot U(\vtG))_{\Nh-k}}=\frac{1}{4}\left(A_{k} + 2iB_{k}\zeta_{2\Nh}^{k} - C_{k}\zeta_{\Nh}^{k}\right)
\end{equation}
From \cref{eqn:modmultABC,eqn:modmultABCconj},
\begin{equation}\label{eqn:modmultSum}
(U(\vtF)\odot U(\vtG))_k + \conj{(U(\vtF)\odot U(\vtG))_{\Nh-k}} = \frac{1}{2}\left(A_k- C_k\zeta_{\Nh}^{k}\right),
\end{equation}
and
\begin{equation}\label{eqn:modmultDiff}
(U(\vtF)\odot U(\vtG))_k - \conj{(U(\vtF)\odot U(\vtG))_{\Nh-k}}= -iB_k\zeta_{2\Nh}^{k}.
\end{equation}
Repacking via \cref{eqn:fourierPack} along with \cref{eqn:modmultSum,eqn:modmultDiff} gives
\begin{align*}
(\vtF\boxdot \vtG)_k 
&= \frac{1}{4}\left(A_k- C_k\zeta_{\Nh}^{k}\right) + \frac{1}{2}\left(B_k\zeta_{2\Nh}^{k}\right)\zeta_{2\Nh}^{-k}= \frac{1}{4}A_k- \frac{1}{4}C_k\zeta_{\Nh}^{k} + \frac{1}{2}B_k.
\end{align*}
Substituting $A_k$, $B_k$, and~$C_k$ and simplifying gives
\begin{align*}
(\vtF\boxdot \vtG)_k &= \frac{1}{4}\left(\tF_k\tG_k + \conj{\tF_{\Nh-k}}\tG_k + \tF_k\conj{\tG_{\Nh-k}}+\conj{\tF_{\Nh-k}}\conj{\tG_{\Nh-k}}\right)-\\
&\ \frac{1}{4}\left(\tF_k\tG_k -\conj{\tF_{\Nh-k}}\tG_k -\tF_k\conj{\tG_{\Nh-k}}+\conj{\tF_{\Nh-k}}\conj{\tG_{\Nh-k}}\right)\zeta_{\Nh}^{k} + \\
&\ \frac{1}{2}\left(\tF_k\tG_k -\conj{\tF_{\Nh-k}}\conj{\tG_{\Nh-k}}\right)\\
&= \tF_k\tG_k - \frac{1}{4}\left(\tF_k - \conj{\tF_{\Nh-k}}\right)\left(\tG_k - \conj{\tG_{\Nh-k}}\right)\left(1+\zeta_{\Nh}^{k}\right),
\end{align*}
as required.

\section*{Acknowledgements}
The authors thank Robert Joseph George for discussions regarding the extension of hybrid dealiasing to arrays of different sizes. The authors also acknowledge the use of ChatGPT and Qwen3.6 during the preparation of this manuscript to assist with code development, correct grammar and typographical errors, and improve the clarity and style of the text. The authors assume responsibility for all content.

Financial support for this work was provided by grants
RES0043585 and RES0046040 from the Natural Sciences and Engineering Research
Council of Canada.

\printbibliography

\end{document}